\documentclass[letterpaper]{cas-dc}
\usepackage[numbers,sort&compress]{natbib}
\usepackage{amsmath,amssymb,amsthm}
\usepackage{graphicx}

\newtheorem{theorem}{Theorem}[section]
\newtheorem{proposition}[theorem]{Proposition}
\newtheorem{lemma}[theorem]{Lemma}
\newtheorem{corollary}[theorem]{Corollary}
\theoremstyle{definition}
\newtheorem{definition}[theorem]{Definition}
\theoremstyle{remark}
\newtheorem{remark}[theorem]{Remark}

\numberwithin{equation}{section}

\newcommand{\R}{\mathbb R}
\newcommand{\N}{\mathcal N}
\shorttitle{Nonlinear KdV--KdV system with time-varying boundary delay}
\shortauthors{R. de A. Capistrano--Filho et al.}

\begin{document}
	\let\WriteBookmarks\relax
	\def\floatpagepagefraction{1}
	\def\textpagefraction{.001}
	
	\title [mode = title]{Exponential boundary stabilization of the KdV--KdV system with time-varying delay: the nonlinear case}
	
	\author[1]{Roberto de A. Capistrano--Filho}[orcid=0000-0002-1608-1832]
	\cormark[1]
	\ead{roberto.capistranofilho@ufpe.br}
	
	\author[2]{Boumedi\`ene Chentouf}[orcid=0000-0001-9365-449X]
	\ead{boumediene.chentouf@ku.edu.kw}
	
	\author[1]{Victor H. Gonzalez~Martinez}[orcid=0000-0002-0090-6908]
	\ead{victor.martinez@ufpe.br}
	
	\author[3]{Juan~Ricardo Mu\~noz}[orcid=0009-0003-9788-792X]
	\ead{juan.munoz@dim.uchile.cl}
	
	\affiliation[1]{organization={Departamento de Matem\'atica, Universidade Federal de Pernambuco},
		addressline={Cidade Universit\'aria},
		postcode={50740-545},
		city={Recife},
		state={PE},
		country={Brazil}}

	\affiliation[2]{organization={Department of Mathematics, College of Science, Kuwait University. Sabah AlSalem University City. Shadadiya. P.O. Box 5969},
		postcode={13060},
		city={Safat},
		country={Kuwait}}
	
	\affiliation[3]{organization={Departamento de Ingenier\'ia Matem\'atica, Facultad de Ciencias F\'isicas y Matem\'aticas, and Centro de Modelamiento Matem\'atico (CMM), Universidad de Chile},
		city={Santiago},
		country={Chile}}
	
	\cortext[cor1]{Corresponding author}
	
	\begin{abstract}
		We consider a Boussinesq system of KdV--KdV type on a bounded interval, with a boundary feedback that includes a damping term with a time-varying delayed term. Although the exponential stability of the linearized system was established in \cite{CCGM} by a Lyapunov approach, global existence and stabilization for the nonlinear problem remained open. Here, we prove that, for sufficiently small initial data and under a suitable restriction on the interval length, the nonlinear system admits a unique global solution whose energy decays exponentially. In particular, we provide an explicit decay rate and a quantitative smallness condition on the initial energy. Our approach combines a Lyapunov argument with well-posedness estimates that are uniform under time translations of the delay, allowing the solutions to be continued globally.
	\end{abstract}
	
	
	\maketitle
	
	\section{Introduction}
	
	The boundary stabilization of the KdV--KdV system with a time-varying delay was studied in \cite{CCGM}, where exponential stability was established for the linearized system and local well-posedness was obtained for the nonlinear problem. In the present note, we improve this analysis by addressing the global existence and exponential stability of the nonlinear solutions for sufficiently small initial data.
	
	Related results on nonlinear boundary stabilization with delay were obtained for the Hirota–Satsuma system in \cite{GMdelay}, where exponential decay for small initial data is established through Lyapunov and observability arguments. A useful structural feature of that system is the cancellation of the nonlinear terms in the natural energy balance, which preserves the dissipativity of the linearized problem. For the original KdV–KdV system, a cubic term of indefinite sign remains in this balance, and the corresponding energy estimate does not directly provide the bound needed for global continuation. In the present note, this difficulty is overcome by retaining the spatial dissipation in the Lyapunov estimate and absorbing the nonlinear contribution for sufficiently small initial energy.
	
	\subsection{Setting}
	Let $L>0$. We study the Boussinesq system of KdV--KdV type
	\begin{equation}\label{eq:KdV}
		\begin{cases}
			\eta_t+\omega_x+\omega_{xxx}+(\eta\omega)_x=0,\\
			\omega_t+\eta_x+\eta_{xxx}+\omega\omega_x=0,\\
			\eta(t,0)=\eta(t,L)=\eta_x(t,0)=0,\\
			\omega(t,0)=\omega(t,L)=0,\\
			\omega_x(t,L)=-\alpha\eta_x(t,L)+\beta\eta_x(t-\tau(t),L),\\
			\eta_x(\theta,L)=h_0(\theta),\qquad \theta\in(-\tau(0),0),\\
			(\eta,\omega)(0,x)=(\eta_0,\omega_0)(x),
		\end{cases}
	\end{equation}
	for $t>0$ and $x\in(0,L)$, which corresponds to system (1.7) in \cite{CCGM}. Here, we recall that $\alpha,\beta\in\mathbb R$ are the feedback gains, $\tau$ is the time-dependent delay, and $h_0$ prescribes
	the delayed state. For the derivation of the KdV--KdV system, we refer to \cite{BCS}; boundary stabilization without delay is studied in \cite{PR,CPR}. Throughout this note, we assume that the delay satisfies
	\begin{equation}\label{eq:tau}
		\begin{cases}
			\tau\in W^{2,\infty}_{\rm loc}([0,\infty)),\qquad 0<\tau_*\leq \tau(t)\le M,\\
			\dot\tau(t)\leq d,\qquad 0\le d<1,
		\end{cases}
	\end{equation}
	for constants $\tau_*,M,d$, and the feedback gains satisfy
	\begin{equation}\label{eq:gain}
		\beta\neq0,\qquad (2\alpha-|\beta|)(1-d)>|\beta|.
	\end{equation}
Note that the second condition in \eqref{eq:gain} is equivalent to \cite[(1.6)]{CCGM}. Moreover, the case $\beta=0$ (no delay) is discussed in Remark \ref{rem:beta0}.
	
	\begin{remark}[On the delay hypothesis]\label{rem:tau}
		The lower bound $\tau(t)\geq\tau(0)$ in \cite[(1.5)]{CCGM} depends on the choice of the initial time. For the translated delay $\tau_s(t)=\tau(s+t)$, the corresponding condition would require $\tau(s+t)\geq\tau(s)$ for $t\geq0$, which does not follow from the original assumption. For example, $\tau(t)=1+\frac14(1-\cos t)$ satisfies $\tau(0)=1\leq\tau(t)\leq\frac32$ and $\dot\tau(t)\leq\frac14$, whereas $\tau_\pi(\pi)=1<\frac32=\tau_\pi(0)$.
		
		In \eqref{eq:tau}, we impose a fixed positive lower bound $\tau_*$, independent of the initial time. Every translated delay $\tau_s$ then satisfies the same bounds with constants $\tau_*,M,d$, and retains the required local regularity. This formulation includes the delays considered in \cite{CCGM}, with $\tau_*=\tau(0)$, and allows the local existence estimates to be applied at each restart with unchanged constants; see Lemma~\ref{lem:shift}.
	\end{remark}
	
	Following \cite{NVF}, we introduce $z(t,\rho)=\eta_x(t-\tau(t)\rho,L)$, $\rho\in(0,1)$, which solves the transport equation
	\begin{gather*}
		\begin{cases}
			\tau(t)z_t(t,\rho)+(1-\dot\tau(t)\rho)z_\rho(t,\rho)=0,\\
			z(t,0)=\eta_x(t,L),\qquad z(0,\rho)=Z_0(\rho),
		\end{cases}
	\end{gather*}
	where $h_0\in L^2(-\tau(0),0)$ is the initial delayed profile of $\eta_x(\cdot,L)$, and
	\begin{equation*}
		Z_0(\rho)=h_0(-\tau(0)\rho),\qquad\rho\in(0,1).
	\end{equation*}
	
	Let us define the initial state by $U_0:=(\eta_0,\omega_0,Z_0)$. From the definition of the delayed profile $Z_0$, we have
	\begin{equation*}
		\tau(0)\|Z_0\|^2_{L^2(0,1)} = \|h_0\|^2_{L^2(-\tau(0),0)}.
	\end{equation*}
	We consider the state space $H=X_0\times L^2(0,1)$, where
	$X_0=[L^2(0,L)]^2$, and define the associated energy by
	\begin{equation}\label{eq:En}
	\begin{split}
		E(t) = \frac12\int_0^L\big(\eta^2+\omega^2\big)(t,x)\,dx\\
		+ \frac{|\beta|}{2}\tau(t)\int_0^1z^2(t,\rho)\,d\rho.
		\end{split}
	\end{equation}
	Consequently, the initial energy can be expressed either in terms
	of $Z_0$ or in terms of the prescribed delayed profile $h_0$
	\begin{align*}
		E(0)&= \frac12\|(\eta_0,\omega_0)\|^2_{X_0} + \frac{|\beta|\tau(0)}{2}\|Z_0\|^2_{L^2(0,1)}\\
		&=\frac12\|(\eta_0,\omega_0)\|^2_{X_0} + \frac{|\beta|}{2}\|h_0\|^2_{L^2(-\tau(0),0)}.
	\end{align*}
	
	When there is no risk of confusion, we write $\|\cdot\|$ for the norm of $L^2(0,L)$ and employ the notation
	\begin{equation*}
		\|U\|^2=\|\eta\|^2+\|\omega\|^2,
		\qquad
		\|\nabla U\|^2=\|\eta_x\|^2+\|\omega_x\|^2.
	\end{equation*}
	Additionally, to have clear computations, let us introduce the boundary quantities
	\begin{gather*}
		a(t):=\eta_x(t,L), \qquad b(t):=z(t,1)=\eta_x(t-\tau(t),L),\\
		v(t)=(a(t),b(t))^T.
	\end{gather*}
	Let us consider
	\begin{equation}\label{eq:Phi_c0}
		\Phi =
		\begin{pmatrix}
			-2\alpha+|\beta| & \beta\\
			\beta & |\beta|(d-1)
		\end{pmatrix},
		\ \
		c_0=-\frac12\lambda_{\max}(\Phi).
	\end{equation}
	Observe that~\eqref{eq:gain} ensures that $\Phi$ is negative definite. Hence $c_0>0$ and is such that
	\begin{equation*}
		\frac12v^T\Phi v\le-c_0|v|^2, \qquad v\in\R^2.
	\end{equation*}
	
Exponential decay of the energy of the linearized system was established in \cite{CCGM} for $0<L<\sqrt3\pi$ through a Lyapunov argument.
	Following that approach, we can consider the same Lyapunov functional
	\begin{equation}\label{eq:V}
		V=E+\mu_1V_1+\mu_2V_2,
	\end{equation}
	where
	\begin{equation}\label{eq:V_i}
		\begin{cases}
			V_1(t)&=\int_0^Lx\,\eta\omega\,dx,\\
			V_2(t)&=\frac{|\beta|}{2}\tau(t)\int_0^1(1-\rho)z^2(t,\rho)\,d\rho.
		\end{cases}
	\end{equation}
	
	For the nonlinear problem, local well-posedness was obtained in \cite[Theorem 2.7]{CCGM}, while global existence and stabilization remained open; see \cite[Remark 2 and Section 4]{CCGM}.
	The main difficulty is to obtain an a priori $L^2$ estimate that accounts for the nonlinear contribution. 
	Indeed, for regular solutions, there holds
	\begin{equation*}
		\frac{d}{dt}\,\frac12\|U\|^2 
		= -\alpha a^2+\beta ab -\frac12\int_0^L\omega_x\eta^2\,dx.
	\end{equation*}
	The nonlinear term has no prescribed sign, and controlling this contribution is the main step in extending the energy estimates to the nonlinear system.
	
	For clarity, we collect in Section \ref{sec:corr} the constants and estimates from \cite{CCGM} used in the present analysis, together with the adjustments needed. Throughout this note, we use the values specified in that section. Let us now introduce the admissibility condition for the parameters of the Lyapunov functional.
	
	\begin{definition}\label{def:adm}
		We say that a pair $(\mu_1,\mu_2)$ is \emph{admissible} if $\mu_1,\mu_2>0$, $\mu_1L<1$, and the matrix
		\begin{equation}\label{eq:Psi}
			\begin{split}
				\Psi = \Phi &+ L\mu_1
				\begin{pmatrix}
					\alpha^2+1 & -\alpha\beta\\
					-\alpha\beta & \beta^2
				\end{pmatrix}+|\beta|\mu_2
				\begin{pmatrix}
					1 & 0\\
					0 & 0
				\end{pmatrix}
			\end{split}
		\end{equation}
		is negative definite. In particular, an explicit family of admissible pairs is given in (C4) of Section \ref{sec:corr}.
	\end{definition}
	
	\subsection{Main results}\label{ssec:main}
	
	Let us now state the main result. Fix $0<L<\sqrt3\pi$ and an admissible pair $(\mu_1,\mu_2)$ in the sense of Definition \ref{def:adm}. We choose the decay parameter $\lambda$ so that
	\begin{equation}\label{eq:lambda}
		0<\lambda< \min\left\{\frac{\mu_1(3\pi^2-L^2)}{L^2(1+\mu_1L)}, \frac{\mu_2(1-d)}{M(1+\mu_2)} \right\},
	\end{equation}
	and introduce the constants
	\begin{equation}\label{eq:cts}
		\begin{cases}
			\delta & = \frac{\mu_1(3\pi^2-L^2)-\lambda L^2(1+\mu_1L)} {2\pi^2}, \\
			K & = \frac14\left(\frac L\pi\right)^{1/2} + \mu_1\left(\frac L\pi\right)^{3/2},\\
			\zeta &= \frac{1+\max\{\mu_1L,\mu_2\}}{1-\mu_1L}.
		\end{cases}
	\end{equation}
The restriction \eqref{eq:lambda} ensures that $\delta>0$. Then, for small enough initial data, we have the following result.
	
	\begin{theorem}[Global existence and exponential stability for small data]
		\label{th:main}
		Assume \eqref{eq:tau}, \eqref{eq:gain}, and $0<L<\sqrt3\pi$, with the parameters chosen as above. 
		If $U_0=(\eta_0,\omega_0,Z_0)\in H$ satisfies
		\begin{equation}\label{eq:smallness}
			E(0)<r^2:=\frac{\delta^2}{2\zeta K^2},
		\end{equation}
		then system \eqref{eq:KdV} admits a unique global solution such that
		\begin{equation*}
			(\eta,\omega) \in C([0,T];X_0) \cap L^2(0,T;H^1_0(0,L)^2)
		\end{equation*}
		for every $T>0$.
		Moreover, its energy satisfies $E(t)\leq\zeta E(0)e^{-\lambda t}$, for $t\ge0$. Finally, setting $\theta= {K\sqrt{2\zeta E(0)}}/{\delta}<1$, we have 
		\begin{equation}\label{eq:mest}
			\begin{aligned}
				&\sup_{t\ge0}\|U(t)\| 
				+ \left(\int_0^\infty\|\nabla U(t)\|^2\,dt \right)^{1/2} \\
				&\quad\leq
				\sqrt{2\zeta E(0)}
				+ \sqrt{\frac{(1+\max\{\mu_1L,\mu_2\})E(0)}{\delta(1-\theta)}}.
			\end{aligned}
		\end{equation}
	\end{theorem}
	
	Theorem \ref{th:main} addresses, for sufficiently small initial data, the question of global existence raised in \cite[Remark 2]{CCGM}. 
	It also extends the linear stabilization result of \cite[Theorem 1.1]{CCGM} to the nonlinear system, with the decay bound specified in \eqref{eq:lambda}.
	The proof combines a Lyapunov estimate with a continuation argument. 
	A key point is that the estimates developed in Section \ref{sec:wp} depend on the delay only through $\tau_*$, $M$, and $d$, the latter through $c_0$.
	Consequently, the same bounds remain available when the initial time is shifted. 
	This allows us to continue the solution in the energy space, without imposing additional regularity on the transport component; see also the question discussed in \cite[Section 4.1(3)]{CCGM}.
	
\begin{figure*}[t]
		\centering
		\includegraphics[width=\textwidth]{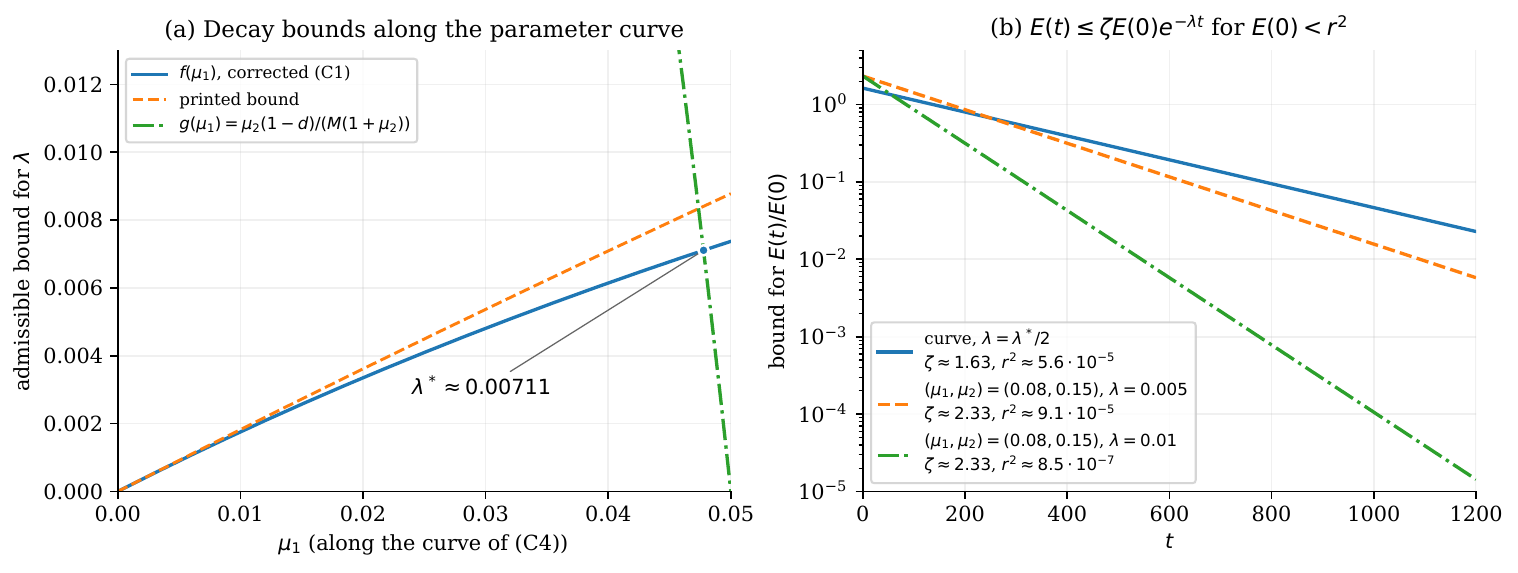}
		\caption{
			Illustration of the decay bounds in Theorem \ref{th:main} for $L=5$, $d=\frac12$, $\alpha=1$, $\beta=\frac12$, and $M=3$.
			(a) Along the parameter curve in (C4), the upper bound for the decay parameter is $\min\{f(\mu_1),g(\mu_1)\}$, where $f$ is given in (C1) and $g(\mu_1)=\mu_2(1-d)/(M(1+\mu_2))$.
			The dashed line corresponds to the expression in \cite[(1.10)]{CCGM}, before the adjustment explained in (C1).
			(b) Bounds $\zeta e^{-\lambda t}$ for $E(t)/E(0)$ associated with three parameter choices, each valid under its corresponding condition $E(0)<r^2$. 
			The comparison illustrates the dependence of the energy bound on the decay rate, the prefactor $\zeta$, and the sufficient smallness threshold.}
		\label{fig:decay}
	\end{figure*}

	To illustrate the dependence of the estimates on the parameters, consider
	\begin{equation*}
		L=5,\quad d=\frac12,\quad \alpha=1,\quad \beta=\frac12,\quad M=3.
	\end{equation*} 
	Along the family of admissible pairs given in (C4), the choice $\mu_1\approx0.0478$, $\mu_2\approx0.0446$, and $\lambda=\lambda^*/2\approx0.00355$ yields
	\begin{equation*}
		\begin{gathered}
			\delta\approx5.58\cdot10^{-3}, \qquad K\approx0.411,\\
			\zeta\approx1.63, \qquad r^2\approx5.65\cdot10^{-5}.
		\end{gathered}
	\end{equation*}
	Here, $\lambda^*$ denotes the maximum of the decay bound along that family; see (C8). 
	On the other hand, the admissible pair $\mu_1=0.08$, $\mu_2=0.15$, with $\lambda=0.005$, gives
	\begin{equation*}
		\begin{gathered}
			\delta\approx9.81\cdot10^{-3}, \qquad K\approx0.476,\\
			\zeta\approx2.33, \qquad r^2\approx9.1\cdot10^{-5}.
		\end{gathered}
	\end{equation*}
	These choices illustrate how both the decay estimate and the sufficient smallness threshold depend on the admissible pair.
	For a fixed pair, increasing $\lambda$ decreases $\delta$ and therefore reduces the radius $r$.
	Figure \ref{fig:decay} illustrates these bounds.
	
	For initial data of arbitrary size, global existence for \eqref{eq:KdV} remains open in the setting considered here.
	The smallness assumption enters through the control of the cubic term $-\frac12\int_0^L\omega_x\eta^2\,dx$, for which the present estimates do not yield a global energy bound for arbitrary data.
	
	\begin{remark}[The case $\beta=0$]\label{rem:beta0}
		When $\beta=0$, the feedback contains no delayed term.
		The state space reduces to $X_0$, and the energy is $E=\frac12\|U\|^2$. Let $\alpha>0$ and choose
		\begin{equation*}
			0<\mu_1<\frac{2\alpha}{L(1+\alpha^2)}.
		\end{equation*}
		Then $\mu_1L<1$ and
		\begin{equation*}
			-\alpha+\frac{L\mu_1}{2}(1+\alpha^2)<0.
		\end{equation*}
		For $0<L<\sqrt3\pi$, we take $V=E+\mu_1V_1$ and choose
		\begin{equation*}
			0<\lambda< \frac{\mu_1(3\pi^2-L^2)}{L^2(1+\mu_1L)}.
		\end{equation*}
		With $\delta$ and $K$ defined as above and
		\begin{equation*}
			\zeta=\frac{1+\mu_1L}{1-\mu_1L},
		\end{equation*}
		the proof of Theorem \ref{th:main} applies after omitting the terms associated with the delay. 
		We therefore obtain the same global existence and exponential decay conclusions for sufficiently small initial data.
	\end{remark}
	
This article is organized as follows: In Section \ref{sec:corr}, we provide some clarifications on the previous work \cite{CCGM} so that the presentation of the current paper becomes clearer and more accurate. In Section \ref{sec:wp}, we prove some estimates for the linearized system that will be useful for the local existence and continuation arguments for the nonlinear problem. Section \ref{nonlin} is devoted to the estimates of the nonlinear term in the problem and then the existence of a unique maximal solution. Section \ref{proofs} is devoted to the proof of the main result Theorem \ref{th:main}. Finally, Section \ref{sec:conclusions} contains some numerical experiments that illustrate Theorem \ref{th:main}, and also concluding remarks.

	\section{Corrections and clarifications to \cite{CCGM}}\label{sec:corr}
	
	Some constants and estimates in \cite{CCGM} contain typographical errors or inconsistencies that require clarification. Since these quantities play an important role in the analysis developed here, we collect below the necessary corrections, together with some refinements concerning the admissibility of the Lyapunov parameters and the optimization of the decay rate. Throughout this note, we systematically use the corrected expressions stated below.
	
	\begin{enumerate}
		
		\item[(C1)] \emph{Decay rate.}
		The estimate of the term $S_1$ in \cite[Section 3]{CCGM} yields the restriction
		\begin{equation*}
			0<\lambda < \frac{\mu_1(3\pi^2-L^2)}{L^2(1+\mu_1L)},
		\end{equation*}
		rather than the expression with denominator $1+\mu_1$ stated in \cite[Theorem 1.1, (1.10)]{CCGM}. In particular, for $L>1$, the printed bound overestimates the decay rate justified by the proof. We observe that the function $f(\mu_1)$ introduced in \cite[Proposition 3.1]{CCGM} already contains the correct denominator.
		
		\item[(C2)] \emph{The functional $V_1$.}
		The calculations in \cite[Section 3]{CCGM} are performed with
		\begin{equation*}
			V_1(t)=\int_0^L x\eta(t,x)\omega(t,x)\,dx,
		\end{equation*}
		without the factor $1/2$ appearing in its original definition. Consequently, the estimate consistent with those calculations is
		\begin{equation*}
			|\mu_1V_1(t)| \leq \frac{\mu_1L}{2}\|U(t)\|^2.
		\end{equation*}
		
		\item[(C3)] \emph{Equivalence of the energy and the Lyapunov functional.}
		Recalling that we are considering the Lyapunov  $V=E+\mu_1V_1+\mu_2V_2$ and $V_2\geq0$, we get
		\begin{equation}\label{eq:EeqV}
			\begin{split}
				(1-\mu_1L)E(t)&\leq V(t)\\
				&\leq \left(1+\max\{\mu_1L,\mu_2\}\right)E(t).
			\end{split}
		\end{equation}
		Thus, under the admissibility condition $\mu_1L<1$, the right constant to be considered is
		\begin{equation*}
			\zeta = \frac{1+\max\{\mu_1L,\mu_2\}}{1-\mu_1L}.
		\end{equation*}
		This improves the constant stated in \cite[(1.9)]{CCGM}, whose denominator is $1-\max\{\mu_1L,\mu_2\}$. In particular, the latter requires the additional restriction $\mu_2<1$ to remain positive, whereas the estimate above holds for every admissible pair. Such a restriction on $\mu_2$ is not guaranteed by the parameter choice described in \cite[Section 4.1]{CCGM}.
		
		\item[(C4)] \emph{Admissible parameters.}
		The argument in \cite[Section 3]{CCGM} requires the matrix $\Psi$, introduced in Definition \ref{def:adm}, to be negative definite. An explicit family of admissible pairs is obtained by choosing
		\begin{equation*}
			0<\mu_1< \frac{(2\alpha-|\beta|)(1-d)-|\beta|}{L(1-d)(1+\alpha^2)}
		\end{equation*}
		and setting
		\begin{equation*}
			\mu_2 = \frac{(2\alpha-|\beta|)(1-d)-|\beta|}{|\beta|(1-d)} -\frac{L(1+\alpha^2)}{|\beta|}\mu_1.
		\end{equation*}
		Indeed, along this curve, a direct computation gives
		\begin{equation*}
			\Psi_{11}=-\frac{|\beta|}{1-d}<0
		\end{equation*}
		and
		\begin{multline*}
			\det\Psi = \frac{L\beta^2\mu_1}{1-d}\left[(1-d)(2\alpha-L\alpha^2\mu_1)-|\beta|\right].
		\end{multline*}
		The gain condition \eqref{eq:gain}, together with the upper bound on $\mu_1$, ensures that the expression in brackets is strictly positive. Hence $\det\Psi>0$, and $\Psi$ is negative definite throughout the prescribed interval.
		
		Moreover, the same choice guarantees $\mu_2>0$ and
		\begin{equation*}
			\mu_1L < \frac{2\alpha}{1+\alpha^2}\le1.
		\end{equation*}
		Therefore, every pair on this curve is admissible, without requiring the additional condition $|\alpha|<1$ imposed in \cite[Section 4.1(5)]{CCGM}. We emphasize, however, that this construction provides only a particular family of admissible pairs and does not characterize the entire admissible set (see (C8) below).
		
		\item[(C5)] \emph{Energy identity}
		For the linearized system, the exact derivative of the energy involves $\dot\tau(t)$ rather than its upper bound $d$. More precisely,
		\begin{equation*}
			\begin{aligned}
				E'(t) &=\frac12v(t)^T\Phi v(t) -\frac{|\beta|}{2}\left(d-\dot\tau(t)\right)b(t)^2\\
				&\leq\frac12v(t)^T\Phi v(t) \leq -c_0\left(a(t)^2+b(t)^2\right).
			\end{aligned}
		\end{equation*}
		Thus, the identity stated in \cite[Proposition 2.3]{CCGM} with the matrix $\Phi$ should be understood as an upper bound. Furthermore, the inequality
		\begin{equation*}
			E'(t)+\eta_x^2(t,L)+z^2(t,1)\leq 0,
		\end{equation*}
		used in the proof of \cite[Proposition 2.4]{CCGM}, must retain the coercivity constant $c_0$, namely,
		\begin{equation*}
			E'(t)+c_0\left(\eta_x^2(t,L)+z^2(t,1)\right)\leq 0.
		\end{equation*}
		
		\item[(C6)] \emph{Weighted a priori estimates.}
		As a consequence of the (C5) correction, the energy and boundary-trace estimates in \cite[Proposition 2.4]{CCGM} must retain the time-dependent weights and the dissipation constant. In particular,
		\begin{multline*}
			\|(\eta,\omega)(t)\|_{X_0}^2 + |\beta|\tau(t)\|z(t)\|_{L^2(0,1)}^2\\
			\leq \|(\eta_0,\omega_0)\|_{X_0}^2 + |\beta|\tau(0)\|Z_0\|_{L^2(0,1)}^2,
		\end{multline*}
		while integration of the corrected dissipation inequality gives
		\begin{align*}
			\int_0^T\big(a(t)^2&+b(t)^2\big)\,dt \leq \frac{E(0)}{c_0}\\
			&=\frac{\|(\eta_0,\omega_0)\|_{X_0}^2+|\beta|\tau(0)\|Z_0\|_{L^2(0,1)}^2}{2c_0}.
		\end{align*}
		
		\item[(C7)] \emph{Linear identity for $V_1$.}
		The term
		\begin{equation*}
			-\int_0^L x\eta\omega\omega_x\,dx
		\end{equation*}
		appearing in \cite[(3.5)]{CCGM} arises from the nonlinear equation and should not be present in the calculation for the linearized system.
		
		\item[(C8)] \emph{On the optimality of $\lambda$.}
		The optimization in \cite[Proposition 3.1]{CCGM} concerns the decay bound restricted to the particular family of admissible pairs described in (C4), rather than the full admissible set. 
		To illustrate this distinction, consider the parameters used in \cite{CCGM},
		\begin{equation*}
			L=5,\quad d=\frac12,\quad \alpha=1,\quad\beta=\frac12,\quad M=3.
		\end{equation*}
		In this case, the curve in (C4) is given by
		\begin{equation*}
			0<\mu_1<\frac1{20}, \qquad \mu_2=1-20\mu_1.
		\end{equation*}
		The maximum of the corrected decay bound along this curve is
		\begin{equation*}
			\lambda^* \approx0.007108903,
		\end{equation*}
		corresponding to $\mu_1\approx0.047772311,$ and $\mu_2\approx0.044553787.$
		
		However, the pair $\mu_1=0.08,$ and $\mu_2=0.15$ is also admissible, since
		\begin{equation*}
			\Psi=
			\begin{pmatrix}
				-0.625&0.3\\
				0.3&-0.15
			\end{pmatrix},
			\quad
			\det\Psi=0.00375>0.
		\end{equation*}
		For this choice, the corrected decay restriction is
		\begin{equation*}
			0<\lambda< \min\{0.010534430,\,0.021739130\},
		\end{equation*}
		which allows, for instance, $\lambda=0.01> \lambda^*$.
		
		In fact, for these parameters,
		\begin{equation*}
			\Psi =
			\begin{pmatrix}
				-\frac32+10\mu_1+\frac12\mu_2&\frac12-\frac52\mu_1\\
				\frac12-\frac52\mu_1& -\frac14+\frac54\mu_1
			\end{pmatrix},
		\end{equation*}
		and
		\begin{equation*}
			\det\Psi = \frac{(5\mu_1-1)(10\mu_1+\mu_2-1)}8.
		\end{equation*}
		Consequently, the full admissible set $\mathcal A$ consists of the pairs $(\mu_1,\mu_2)$ with
		\begin{equation*}
			0<\mu_1<\frac1{10},\qquad 0<\mu_2<1-10\mu_1.
		\end{equation*}
		Since the decay bound is nondecreasing with respect to $\mu_2$, its supremum over
		$\mathcal A$ is
		\begin{equation*}
			\begin{split}
				\Lambda_{\mathrm{adm}}
				&=\max_{0\leq s\leq1/10}
				\min\left\{ \frac{s(3\pi^2-25)}{25(1+5s)}, \frac{1-10s}{6(2-10s)}\right\}\\
				&\approx0.011658099.
			\end{split}
		\end{equation*}
		The first function is strictly increasing, whereas the second is strictly decreasing.
		Their unique intersection occurs at $s_*\approx0.092479062$, and determines the above supremum. Since
		the corresponding limiting parameter $\mu_2=1-10s_*$ satisfies $\det\Psi=0$, the supremum is not attained by any admissible pair.
		
		Therefore, the optimization in \cite[Proposition 3.1]{CCGM} yields the largest decay bound along the prescribed parameter curve, but not over the full admissible set. Moreover, this optimization concerns the decay estimates obtained by the Lyapunov method and does not establish an optimal exponential decay rate for the underlying evolution system.
		
	\end{enumerate}
	
	With these corrections, the exponential-decay conclusion of \cite[Theorem 1.1]{CCGM} remains valid for the linearized system, with the constants and admissibility conditions specified above. In the next section, we revisit the well-posedness estimates and establish bounds that are uniform under time translations of the delay, which will be essential for the analysis of the nonlinear problem.
	
	\section{Well-posedness with constants uniform under time translations}
	\label{sec:wp}
	
	In this section, we establish the estimates for the linearized system that will be used in the local existence and continuation arguments for the nonlinear problem. First, let us set the time-dependent norm such that for each $t\geq 0$, we equip $H$ with the norm
	\begin{equation*}
		\|U\|_t^2 = \|\eta\|^2+\|\omega\|^2 + |\beta|\tau(t)\|z\|_{L^2(0,1)}^2,
	\end{equation*}
	for $U=(\eta,\omega,z)$.
	Thus, along a solution $U$, we have: $\|U(t)\|_t^2=2E(t)$. Note that, by setting
	\begin{equation*}
		c_*=\min\{1,|\beta|\tau_*\}, \qquad C_*=\max\{1,|\beta|M\},
	\end{equation*}
	we get
	\begin{equation*}
		c_*\|U\|_H^2 \leq \|U\|_t^2\leq C_*\|U\|_H^2,
	\end{equation*}
	for all $U\in H$ and $t\geq 0$.
	
	For $0\leq s < t_1$, let us set the solution spaces
	\begin{equation*}
		B_{s,t_1} =C([s,t_1];X_0) \cap L^2\left(s,t_1;H_0^1(0,L)^2\right),
	\end{equation*}
	endowed with the norm
	\begin{multline*}
		\|(\eta,\omega)\|_{B_{s,t_1}}
		= \sup_{t\in[s,t_1]}\|U(t)\| + \left(\int_s^{t_1}\|\nabla U(t)\|^2\,dt\right)^{1/2}.
	\end{multline*}
	Similarly, we write $B_T=B_{0,T}$ and use the same notation for scalar functions, that is,
	\begin{equation*}
		\|u\|_{B_T} = \sup_{t\in[0,T]}\|u(t)\| + \|u_x\|_{L^2(0,T;L^2(0,L))}.
	\end{equation*}

		Recall that the linearized system is governed by the family $\{A(t)\}$ defined by (see \cite[(2.4)--(2.5)]{CCGM})
		\begin{multline*}
		A(t)\left(\eta, \omega , z\right) :=\Big(-\omega_x-\omega_{xxx}, -\eta_x-\eta_{xxx},\\
		 \frac{\dot\tau(t)\rho -1}{\tau(t)} z_\rho\Big),
		\end{multline*}
		whose domain $D(A(t))$ consists of the triples
		\begin{equation*}
			(\eta,\omega,z) \in \left[ H^3(0,L)\cap H_0^1(0,L) \right]^2 \times H^1(0,1)
		\end{equation*}
		such that
		\begin{gather*}
			\eta_x(0)=0,\qquad z(0) = \eta_x(L),\\
			\omega_x(L) = -\alpha \eta_x(L) + \beta z(1).
		\end{gather*}
		Clearly, $D(A(t))=D(A(0))$ is independent of $t$, and Kato's theory \cite{Kato70} applies under the conditions \eqref{eq:tau}--\eqref{eq:gain}. Hence $\{A(t)\}$ generates an evolution system $\{S(t,s)\}_{0\leq s\leq t}$ on $H$; see \cite[Theorem 2.2]{CCGM}.

	Now, aiming at the nonlinear problem, we set the mild solution for a system with a source. Let $s\geq 0$, given $U_s\in H$ and $f=(f_1,f_2)\in L^1(s,t_1;X_0)$, we consider
	\begin{equation}\label{eq:lin.sou}
		\begin{cases}
			\eta_t+\omega_x+\omega_{xxx}=f_1,
			\\
			\omega_t+\eta_x+\eta_{xxx}=f_2,
		\end{cases}
	\end{equation}
	coupled with the transport equation for $z$ and the boundary conditions of \eqref{eq:KdV}, with initial state $U_s$ at time $s$. Its mild solution is then given by the Duhamel formula
	\begin{equation*}
		U(t)=S(t,s)U_s + \int_s^t S(t,r)\left(f_1(r),f_2(r),0\right)\,dr,
	\end{equation*}
	for $t\in[s,t_1]$.
	
	\begin{proposition}[Linear estimates]\label{pr:lin}
		Assume \eqref{eq:tau}--\eqref{eq:gain}. Let $s\ge0$, $T>0$, $U_s\in H$, and $f\in L^1(s,s+T;X_0)$, and set $F:=\|f\|_{L^1(s,s+T;X_0)}$. Then the mild solution on $[s,s+T]$ satisfies
		$(\eta,\omega)\in B_{s,s+T}$ and
	\begin{gather*}
				z\in C([s,s+T];L^2(0,1)),\ \ 
				a,b\in L^2(s,s+T).
		\end{gather*}
		Moreover, we have
		\begin{align*}
			& \sup_{t\in[s,s+T]}\|U(t)\|_t \leq\|U_s\|_s+2F,\\
			& c_0\int_s^{s+T}(a^2+b^2)\,dt \leq\|U_s\|_s^2+3F^2,\\
			& \int_s^{s+T}\|\nabla U(t)\|^2\,dt \leq \Lambda(T) \left(\|U_s\|_s+F\right)^2,
		\end{align*}
		where
		\begin{equation*}
			\Lambda(T) = \frac23
			\left[ 2T+6L+\frac{3L(1+2\alpha^2+2\beta^2)}{2c_0}\right].
		\end{equation*}
		In particular, with $C_T=2+\sqrt{\Lambda(T)}$,
		\begin{equation*}
			\|(\eta,\omega)\|_{B_{s,s+T}} \leq C_T\left(\|U_s\|_s+F\right).
		\end{equation*}
		The constant $C_T$ is nondecreasing in $T$, independent of $s$, and uniform over the delays satisfying \eqref{eq:tau}. In addition, for every $t\in[s,s+T]$, 
		\begin{equation}\label{eq:Id1}
			\begin{aligned}
				&\frac12\|U(t)\|_t^2-\frac12\|U_s\|_s^2\\
				&\ =\int_s^t\Big[\Big(-\alpha+\frac{|\beta|}{2}\Big)a^2+\beta ab-\frac{|\beta|}{2}(1-\dot\tau)b^2\Big]dr\\
				&\quad+\int_s^t\int_0^L(\eta f_1+\omega f_2)\,dx\,dr,
			\end{aligned}
		\end{equation}
		\begin{equation}\label{eq:Id2}
			\begin{aligned}
				&V_1(t)-V_1(s)
				=\int_s^t\Big[\frac12\|U\|^2-\frac32\|\nabla U\|^2\Big]dr\\
				&\quad+\frac L2\int_s^t\left((\alpha a-\beta b)^2+a^2\right)dr\\
				&\quad+\int_s^t\int_0^L x(\omega f_1+\eta f_2)\,dx\,dr,
			\end{aligned}
		\end{equation}
		and
		\begin{equation}\label{eq:Id3}
			\begin{split}
				&V_2(t)-V_2(s)\\
				&\ = \frac{|\beta|}{2}\int_s^t\Big[ a(r)^2-\int_0^1(1-\dot\tau(r)\rho)z(r,\rho)^2\,d\rho \Big]dr.
			\end{split}
		\end{equation}
	\end{proposition}
	
	\begin{proof}
		Let us split the proof in several steps. 
		
		\medskip
		\noindent\textit{Step 1. Identities for regular data.}
		Let $U_s\in D(A(0))$ and $f\in C^1([s,s+T];C_c^\infty(0,L)^2)$. Then,  $(f_1,f_2,0)$ belongs to $C^1([s,s+T];D(A(0)))$, and consequently the corresponding solution is classical; see \cite{Kato70,Pazy}, that is,
		\begin{equation*}
			U\in C([s,s+T];D(A(0)))\cap C^1([s,s+T];H).
		\end{equation*}
		
		Multiplying the first equation by $\eta$ and the second by $\omega$, integrating over $(0,L)$, and using the boundary conditions we obtain
		\begin{equation}\label{eq:K1}
			\frac{d}{dt}\frac12\|U\|^2 = -\alpha a^2+\beta ab + \int_0^L(\eta f_1+\omega f_2)\,dx.
		\end{equation}
		
		On the other hand, the transport equation yields
		\begin{equation}
			\begin{aligned}\label{eq:K2}
				&\frac{d}{dt}\left( \frac{|\beta|}{2}\tau\int_0^1z^2\,d\rho \right)\\
				&\ =\frac{|\beta|}{2}\left[\dot\tau\int_0^1z^2\,d\rho - \int_0^1(1-\dot\tau\rho)(z^2)_\rho\,d\rho \right]\\
				&\ =\frac{|\beta|}{2}\left(a^2-(1-\dot\tau)b^2\right).
			\end{aligned}
		\end{equation}
		Adding~\eqref{eq:K1} and~\eqref{eq:K2} and integrating from $s$ to $t$ proves \eqref{eq:Id1}. For $V_2$, the same calculation with the weight $1-\rho$ gives
		\begin{align*}
			V_2'
			&=\frac{|\beta|}{2}
			\Big[\dot\tau\int_0^1(1-\rho)z^2\,d\rho\\
			&\qquad-\int_0^1(1-\rho)(1-\dot\tau\rho)(z^2)_\rho\,d\rho
			\Big]\\
			&=\frac{|\beta|}{2}\left[a^2-\int_0^1(1-\dot\tau\rho)z^2\,d\rho
			\right].
		\end{align*}
		Integrating from $s$ to $t$ proves \eqref{eq:Id3}.
		
		Next, multiplying the first equation by $x\omega$ and the second by $x\eta$, integration by parts implies
		\begin{align*}
			\int_0^L x\eta_t\omega\,dx
			&=\frac12\|\omega\|^2-\frac32\|\omega_x\|^2 +\frac L2\omega_x(t,L)^2\\
			&\quad+\int_0^L xf_1\omega\,dx,\\
			\int_0^L x\eta\omega_t\,dx
			&=\frac12\|\eta\|^2-\frac32\|\eta_x\|^2 + \frac L2a^2\\
			&\quad+\int_0^L x\eta f_2\,dx.
		\end{align*}
		Adding the results and utilizing the boundary condition implies~\eqref{eq:Id2}. 
		
		\medskip
		\noindent\textit{Step 2. Uniform estimates.}
		By (C5) and~\eqref{eq:Phi_c0}, the boundary terms in \eqref{eq:Id1} satisfy
		\begin{multline*}
			\left(-\alpha+\frac{|\beta|}{2}\right)a^2 + \beta ab-\frac{|\beta|}{2}(1-\dot\tau)b^2\leq -c_0(a^2+b^2).
		\end{multline*}
		Set $Y:=\sup_{r\in[s,s+T]}\|U(r)\|_r$. Since $\|U(t)\|\le\|U(t)\|_t$, the Cauchy--Schwarz inequality and \eqref{eq:Id1} imply
		\begin{equation*}
			\frac12\|U(t)\|_t^2 +c_0\int_s^t(a^2+b^2)\,dr \leq \frac12\|U_s\|_s^2 + YF,
		\end{equation*}
		for every $t\in[s,s+T]$. Dropping the integral and taking the supremum in $t$, we get $Y^2-2FY-\|U_s\|_s^2\leq0$, and hence
		\begin{equation*}
			Y\leq F + \sqrt{F^2+\|U_s\|_s^2}\leq \|U_s\|_s+2F.
		\end{equation*}
		Taking $t=s+T$ in the same inequality yields
		\begin{align*}
			c_0\int_s^{s+T}(a^2+b^2)\,dt
			&\leq\frac12\|U_s\|_s^2 +YF\\
			&\leq\frac12\|U_s\|_s^2 + \left(\|U_s\|_s+2F\right) F\\
			&\leq\|U_s\|_s^2+3F^2.
		\end{align*}
		
		From \eqref{eq:Id2} together with
		\begin{gather*}
			|V_1(t)|\leq\frac L2\|U(t)\|^2,
			\qquad
			(\alpha a-\beta b)^2\le2\alpha^2a^2+2\beta^2b^2,\\
			\left|\int_0^L x(\omega f_1+\eta f_2)\,dx\right| 
			\leq L\|U\|\,\|f\|_{X_0},
		\end{gather*}
		we obtain
		\begin{align*}
			&\frac32\int_s^{s+T}\|\nabla U\|^2\,dt
			\leq\frac T2Y^2+LY^2+LYF\\
			&\qquad+\frac L2(1+2\alpha^2+2\beta^2)\int_s^{s+T}(a^2+b^2)\,dt.
		\end{align*}
		Since $Y\leq2\left(\|U_s\|_s+F\right)$ and
		\begin{equation*}
			\|U_s\|_s^2+3F^2\leq3\left(\|U_s\|_s+F\right)^2,
		\end{equation*}
		we obtain
		\begin{equation*}
			\frac32\int_s^{s+T}\|\nabla U\|^2\,dt
			\leq\frac32\,\Lambda(T)\left(\|U_s\|_s+F\right)^2.
		\end{equation*}
		This proves the estimate with $\Lambda(T)$, and the bound in $B_{s,s+T}$ follows immediately. Notice that the estimates are uniform, i.e., none of these constants depends on $s$ or on the local derivatives of $\tau$ beyond the fixed upper bound $d$ for $\dot\tau$.
		
		\medskip
		\noindent\textit{Step 3. Passage to general data.}
		Let $U_s\in H$ and $f\in L^1(s,s+T;X_0)$ be arbitrary. Since $D(A(0))$ is dense in $H$ and $C^1([s,s+T];C_c^\infty(0,L)^2)$ is dense in $L^1(s,s+T;X_0)$, there exist $U_s^n\in D(A(0))$ and $f^n\in C^1([s,s+T];C_c^\infty(0,L)^2)$ such that
		\begin{equation*}
			U_s^n\to U_s\ \text{in }H,
			\qquad
			f^n\to f\ \text{in }L^1(s,s+T;X_0).
		\end{equation*}
		Let $U^n=(\eta^n,\omega^n,z^n)$ be the corresponding classical solutions, with traces $a^n,b^n$. By linearity, $U^n-U^m$ is the classical solution with initial state $U_s^n-U_s^m$ and source $f^n-f^m$. Hence the estimates of \textit{Step 2}, applied to $U^n-U^m$, together with the uniform equivalence of the norms, show that $(U^n)$, $(\eta^n,\omega^n)$ and $(a^n,b^n)$ are Cauchy sequences, so that
		\begin{align*}
			U^n\to U& \quad\text{in }C([s,s+T];H),\\
			(\eta^n,\omega^n)\to(\eta,\omega) & \quad\text{in }B_{s,s+T},\\
			(a^n,b^n)\to(a,b)& \quad\text{in }L^2(s,s+T;\mathbb R^2).
		\end{align*}
		
		We claim that $U$ is the mild solution with initial state $U_s$ and source $f$. Indeed, each $U^n$ satisfies the Duhamel formula with data $(U_s^n,f^n)$, and \textit{Step 2} with $f=0$ gives $\|S(t,r)\|_{\mathcal L(H)}\leq(C_*/c_*)^{1/2}$ for $s\leq r\leq t\leq s+T$; thus both sides of the formula converge in $C([s,s+T];H)$. Moreover, $(\eta,\omega)(t)\in H_0^1(0,L)^2$ for almost every $t$, since $L^2(s,s+T;H_0^1(0,L)^2)$ is complete.
		
		Finally, the identities \eqref{eq:Id1}--\eqref{eq:Id3} hold for each $U^n$, and every term in them passes to the limit: the terms that are quadratic in $(a,b)$ and in $\nabla U$ by the $L^2$ convergences above; the terms involving $z$ because $z^n\to z$ in $C([s,s+T];L^2(0,1))$ and $\dot\tau\in L^\infty(s,s+T)$; and the source terms because $U^n\to U$ uniformly in time while $f^n\to f$ in $L^1(s,s+T;X_0)$. Therefore \eqref{eq:Id1}--\eqref{eq:Id3} hold for $U$ and every $t\in[s,s+T]$, and the estimates of \textit{Step 2}, which only use these identities, remain valid for $U$.
	\end{proof}
	
It is noteworthy that the preceding estimates also give a uniform bound in the fixed norm
	of $H$.
	\begin{corollary}[Uniform estimate]\label{cor:unif}
		Under the assumptions of Proposition \ref{pr:lin}, we have
		\begin{align*}
			&\sup_{t\in[s,s+T]}\|U(t)\|_H +\|(a,b)\|_{L^2(s,s+T;\mathbb R^2)}\\
			&\quad+\left(\int_s^{s+T} \left(\|U(t)\|^2+\|\nabla U(t)\|^2\right)dt\right)^{1/2}\\
			&\qquad\leq 
			\widetilde C_T\left(\|U_s\|_H+F\right),
		\end{align*}
		where
		\begin{equation*}
			\begin{split}
				\widetilde C_T = \sqrt{\frac{C_*}{c_*}}
				\Bigg[ 2&+\sqrt{\left(1+\frac{L^2}{\pi^2}\right)\Lambda(T)}+\sqrt{\frac3{c_0}} \Bigg].
			\end{split}
		\end{equation*}
		This constant depends only on $T,L,\alpha,\beta,d,\tau_*,M$ and is independent of $s$.
	\end{corollary}
	
	\begin{proof}
		The equivalence of the norms gives
		\begin{equation*}
			\|U(t)\|_H\leq c_*^{-1/2}\|U(t)\|_t,
			\qquad
			\|U_s\|_s\leq C_*^{1/2}\|U_s\|_H.
		\end{equation*}
		Moreover, Poincar\'e's inequality yields
		\begin{equation*}
			\|U(t)\|^2\le\frac{L^2}{\pi^2}\|\nabla U(t)\|^2,
		\end{equation*}
		for a.e.\ $t\in(s,s+T)$.
		Applying Proposition \ref{pr:lin}, we therefore bound the left-hand
		side by
		\begin{multline*}
			\left[\frac2{\sqrt{c_*}} +\sqrt{\left(1+\frac{L^2}{\pi^2}\right)\Lambda(T)} +\sqrt{\frac3{c_0}}\right]\left(\|U_s\|_s+F\right).
		\end{multline*}
		The conclusion follows from $c_*\leq1\leq C_*$ and
		\begin{equation*}
			\|U_s\|_s+F
			\leq C_*^{1/2} \left(\|U_s\|_H+F\right).
		\end{equation*}
	\end{proof}
	
	\begin{remark}[Traces in the energy space]\label{rem:traces}
		Notice that $(\eta,\omega)\in B_T$ alone does not provide a classical trace for $\eta_x(\cdot,L)$. Indeed, the traces $a$ and $b=z(\cdot,1)$ in Proposition \ref{pr:lin} are defined as $L^2(0,T)$ limits of regular boundary traces. The linear estimates ensure that these limits are independent of the approximation and agree with the classical traces
		for regular solutions.
		
		Let $\widehat a$ denote the extension of $a$ to $(-\tau(0),T)$ by
		the initial delayed profile $h_0(\sigma)=Z_0(-\sigma/\tau(0))$ on $(-\tau(0),0)$. Then
		\begin{equation*}
			z(t,\cdot)=\widehat a(t-\tau(t)\,\cdot)\quad\text{in }L^2(0,1),
		\end{equation*}
		for $t\in[0,T]$, and $b(t)=\widehat a(t-\tau(t))$ for almost every $t\in(0,T)$.
		Indeed, $\tau(t)\ge\tau_*$ and $1-\dot\tau(t)\ge1-d>0$ give,
		by changes of variables,
		\begin{gather*}
			\sup_{t\in[0,T]}
			\int_0^1|\widehat a(t-\tau(t)\rho)|^2\,d\rho
			\leq\frac1{\tau_*}\|\widehat a\|_{L^2(-\tau(0),T)}^2,\\
			\int_0^T|\widehat a(t-\tau(t))|^2\,dt
			\leq\frac1{1-d}\|\widehat a\|_{L^2(-\tau(0),T)}^2.
		\end{gather*}
		Applied to differences of regular approximations, these bounds extend the characteristic formulas to all initial states in $H$, uniformly in time in $L^2(0,1)$ for $z$ and in $L^2(0,T)$ for $b$.
	\end{remark}
	
	We conclude with the time-translation property that will allow us to
	restart the nonlinear problem from any state reached by its evolution.
	
	\begin{lemma}[Time translations]\label{lem:shift}
		Let $U$ be a mild solution of the linear problem with source $f\in L^1(0,t_1;X_0)$ on $[0,t_1]$, and let $0\leq s<t_1$. Then $\tau_s(t)=\tau(s+t)$ satisfies \eqref{eq:tau} with the same constants $\tau_*,M,d$. Moreover, $t\mapsto U(s+t)$ is the mild solution on $[0,t_1-s]$ corresponding to the delay $\tau_s$, the source $f(s+\cdot)$, and the initial state $U(s)\in H$.
	\end{lemma}
	
	\begin{proof}
		The assumptions in \eqref{eq:tau} are preserved by time translation. Since $\dot\tau\le d<1$, we have $t-\tau(t)\ge s-\tau(s)$ for $t\ge s$. By Remark \ref{rem:traces},
		$z(s,\rho)=\widehat a(s-\tau(s)\rho)$ for a.e.\ $\rho\in(0,1)$, so $z(s)$ becomes the initial delayed profile to restart the problem at time $s$.
		
		The operator family associated with $\tau_s$ is $\{A(s+t)\}_{t\geq0}$, with domain $D(A(0))$ and evolution operators $S(s+t,s+r)$. Splitting the Duhamel integral at $s$, we get
		\begin{multline*}
			U(s+t)=S(s+t,s)U(s)\\
			+ \int_0^t S(s+t,s+r) \left (f_1(s+r),f_2(s+r),0\right)\,dr,
		\end{multline*}
		for $0\leq t\leq t_1-s$.
		This is the mild formulation of the problem with delay $\tau_s$, source $f(s+\cdot)$, and initial state $U(s)$.
	\end{proof}
	
	\section{The nonlinear system}\label{nonlin}
	
	In order to consider the nonlinear terms, we first estimate the terms in the source space of Proposition \ref{pr:lin}.
	
	\begin{lemma}\label{lem:bil}
		For every $w\in H_0^1(0,L)$,
		\begin{equation}\label{sob}
			\|w\|_{L^\infty}^2\le\|w\|\,\|w_x\|.
		\end{equation}
		Moreover, for scalar $u,v\in B_T$,
		\begin{equation}\label{eq:iq.prod}
			\int_0^T\|(uv)_x(t)\|\,dt \leq 2T^{1/4}\|u\|_{B_T}\|v\|_{B_T}.
		\end{equation}
		Consequently, the nonlinear operator
		\begin{equation*}
			\N(\eta,\omega)=\left((\eta\omega)_x,\omega\omega_x\right)
		\end{equation*}
		satisfies, for $(\eta_i,\omega_i)\in B_T$, $i=1,2$,
		\begin{align*}
			&\|\N(\eta_1,\omega_1)-\N(\eta_2,\omega_2)\|_{L^1(0,T;X_0)}\\
			&\quad\le3T^{1/4} \left(\|(\eta_1,\omega_1)\|_{B_T} + \|(\eta_2,\omega_2)\|_{B_T}\right)\\
			&\qquad\times \|(\eta_1-\eta_2,\omega_1-\omega_2)\|_{B_T},
		\end{align*}
		and
		\begin{equation}\label{eq:NL.est}
			\|\N(\eta,\omega)\|_{L^1(0,T;X_0)} 
			\leq 3T^{1/4}\|(\eta,\omega)\|_{B_T}^2.
		\end{equation}
		The same holds on every interval $[s,s+T]$.
	\end{lemma}
	
	\begin{proof}
		Since $w(0)=w(L)=0$,
		\begin{align*}
			w^2(x)
			&=\int_0^x w(y)w'(y)\,dy-\int_x^L w(y)w'(y)\,dy\\ &\leq \|w\|\,\|w_x\|.
		\end{align*}
		Thus, inequality \eqref{sob} gives
		\begin{equation*}
			\|(uv)_x\|
			\leq\|u_x\|\,\|v\|^{\frac12}\|v_x\|^{\frac12} + \|u\|^{\frac12}\|u_x\|^{\frac12}\|v_x\|.
		\end{equation*}
		Then, applying H\"older's inequality after integrating in time,
		\begin{align*}
			&\int_0^T\|u_x(t)\|\,\|v(t)\|^{1/2}\|v_x(t)\|^{1/2}\,dt\\
			&\quad
			\leq T^{\frac14}
			\|u_x\|_{L^2(0,T;L^2)} \|v\|_{C([0,T];L^2)}^{\frac12} \|v_x\|_{L^2(0,T;L^2)}^{\frac12}\\
			&\quad
			\leq T^{1/4}\|u\|_{B_T}\|v\|_{B_T},
		\end{align*}
		where $L^2=L^2(0,L)$.
		Exchanging $u$ and $v$ gives~\eqref{eq:iq.prod}. On the other hand, for the nonlinear operator $\N$, we can rewrite the terms as
		\begin{gather*}
		\begin{cases}
			\eta_1\omega_1-\eta_2\omega_2
			=\eta_1(\omega_1-\omega_2)+(\eta_1-\eta_2)\omega_2,\\
			\omega_1\omega_{1,x}-\omega_2\omega_{2,x}
			=\frac12\big((\omega_1+\omega_2)(\omega_1-\omega_2)\big)_x,
			\end{cases}
		\end{gather*}
		and applying~\eqref{eq:iq.prod} to each term implies~\eqref{eq:NL.est}. Notice that time translation gives the same bounds on $[s,s+T]$.
	\end{proof}
	
	A solution of \eqref{eq:KdV} on $[0,T]$ is a state $U=(\eta,\omega,z)\in C([0,T];H)$ with $(\eta,\omega)\in B_T$ that satisfies Duhamel's formula of the linear problem \eqref{eq:lin.sou} with initial state $U_0$ and source $-\N(\eta,\omega)\in L^1(0,T;X_0)$.
	
	\begin{theorem}\label{th:local}
		Assume \eqref{eq:tau}--\eqref{eq:gain} and $L>0$.
		\begin{enumerate}
			\item[(a)] Let $T>0$. If $24C_T^2T^{1/4}\|U_0\|_0\leq 1$, then system \eqref{eq:KdV} has a unique solution on $[0,T]$, satisfying
			\begin{equation}\label{eq:nl.s.bd}
				\|(\eta,\omega)\|_{B_T}\leq 2C_T\|U_0\|_0.
			\end{equation}
			The data-to-solution map is Lipschitz on this set of initial states, with values in
			\begin{equation*}
				B_T\times C([0,T];L^2(0,1)).
			\end{equation*}
			
			\item[(b)] For every $U_0\in H$, the conclusion of (a) holds on
			$[0,T_*(\|U_0\|_0)]$, where
			\begin{equation*}
				T_*(\|U_0\|_0)=
				\begin{cases}
					\min\left\{1,(24C_1^2\|U_0\|_0)^{-4}\right\},\\
					\qquad\qquad\qquad\text{if }\|U_0\|_0>0,\\
					1,\qquad\qquad\ \ \text{if }\|U_0\|_0=0.
				\end{cases}
			\end{equation*}
			
			\item[(c)] There is a unique maximal solution on $[0,T_{\max})$. If $T_{\max}<\infty$, then
			\begin{equation*}
				\sup_{0\leq t<T_{\max}}\|U(t)\|_t=\infty.
			\end{equation*}
		\end{enumerate}
	\end{theorem}
	
	We highlight that the same conclusions hold for the problem started at any time $s\geq 0$ from $U_s\in H$, with $\|U_0\|_0$ replaced by $\|U_s\|_s$ and the same constants.
	
	\begin{proof}
		{(a)} Let us consider $R = 2C_T\|U_0\|_0$ be the radius of the closed ball 
		\begin{equation*}
			\mathcal B_R = \left\{ (\eta,\omega) \in B_T \colon \|(\eta,\omega)\|_{B_T} \leq R \right\}.
		\end{equation*}
		
		For $(\eta,\omega)\in \mathcal B_R$, let $\Gamma(\eta,\omega)$ denote the first two components of the mild linear solution with initial state $U_0$ and source $-\N(\eta,\omega)$. By Proposition~\ref{pr:lin}, Lemma~\ref{lem:bil}, and the assumption in (a),
		\begin{align*}
			\|\Gamma(\eta,\omega)\|_{B_T}
			&\leq C_T\|U_0\|_0 + 3C_TT^{1/4}\|(\eta,\omega)\|_{B_T}^2\\
			&\leq C_T\|U_0\|_0 + 12C_T^3T^{1/4}\|U_0\|_0^2\\
			&\leq 2C_T\|U_0\|_0,
		\end{align*}
		and
		\begin{align*}
			&\|\Gamma(\eta_1,\omega_1)-\Gamma(\eta_2,\omega_2)\|_{B_T}\\
			&\quad\leq 12C_T^2T^{1/4}\|U_0\|_0\|(\eta_1-\eta_2,\omega_1-\omega_2)\|_{B_T}\\
			&\quad\leq \frac12\|(\eta_1-\eta_2,\omega_1-\omega_2)\|_{B_T}.
		\end{align*}
		Thus $\Gamma$ is a contraction of the ball into itself, and Banach's fixed point theorem yields a mild solution satisfying~\eqref{eq:nl.s.bd}.
		
		To prove uniqueness, let $U=(\eta,\omega,z)$ and $\widetilde U=(\widetilde\eta,\widetilde\omega,\widetilde z)$ be two mild solutions on $[0,T]$ with the same initial state. For $0\leq s<t\leq T$ with $U(s)=\widetilde U(s)$, Proposition~\ref{pr:lin} and Lemmas~\ref{lem:bil}
		and~\ref{lem:shift} give 
		\begin{align*}
			&\|(\eta-\widetilde\eta,\omega-\widetilde\omega)\|_{B_{s,t}}\\
			&\ \leq 3C_T(t-s)^{\frac14}
			\left(\|(\eta,\omega)\|_{B_T} +\|(\widetilde\eta,\widetilde\omega)\|_{B_T}\right)\\
			&\quad\times
			\|(\eta-\widetilde\eta,\omega-\widetilde\omega)\|_{B_{s,t}}.
		\end{align*}
		For $t-s$ sufficiently small, one can claim that  $U=\widetilde U$ on $[s,t]$.
		Starting at $s=0$ and iterating over finitely many such intervals
		proves uniqueness on $[0,T]$.
		
		Now, given two initial states satisfying (a), their difference can be estimated as follows
		\begin{multline*}
			\|(\eta-\widetilde\eta,\omega-\widetilde\omega)\|_{B_T}
			\leq C_T\|U_0-\widetilde U_0\|_0
			+ \frac12\|(\eta-\widetilde\eta,\omega-\widetilde\omega)\|_{B_T},
		\end{multline*}
		and hence Lipschitz property follows.
		
		\medskip
		\noindent {(b)} Since $C_T\leq C_1$ for $T\leq1$, the stated choice of $T_*(\|U_0\|_0)$ satisfies the condition in (a). For $U_0=0$, the solution is identically zero. The conclusions for a starting time $s$ follow from Lemma \ref{lem:shift} and the uniformity of the linear estimates.
		
		\medskip
		\noindent{(c)} Uniqueness allows the local solutions to be joined into a maximal solution. Suppose that $T_{\max}<\infty$ and
		\begin{equation*}
			\sup_{0\leq t<T_{\max}}\|U(t)\|_t<\infty.
		\end{equation*}
		By (b), every restart from $U(s)$, $0\leq s<T_{\max}$, has an existence time bounded below by
		\begin{equation*}
			T_*\left(\sup_{0\leq t<T_{\max}}\|U(t)\|_t\right)>0.
		\end{equation*}
		Choosing $s<T_{\max}$ sufficiently close to $T_{\max}$, one extends the solution beyond $T_{\max}$, a contradiction.
	\end{proof}
	
	\section{Proof of Theorem \ref{th:main}} \label{proofs}
	
	Here, we are going to estimate the isolated identities that are involved in the Lyapunov estimates. Recall that the linear ones are already estimated.
	\begin{lemma}\label{lem:NL}
		For $\eta,\omega\in H_0^1(0,L)$ and $f=-\N(\eta,\omega)$,
		\begin{gather*}
			\int_0^L(\eta f_1+\omega f_2)\,dx = -\frac12\int_0^L\omega_x\eta^2\,dx,\\
			\int_0^Lx(\omega f_1+\eta f_2)\,dx = \int_0^L\eta\omega^2\,dx,
		\end{gather*}
		and
		\begin{equation*}
			\left|-\frac12\int_0^L\omega_x\eta^2\,dx +\mu_1\int_0^L\eta\omega^2\,dx\right| \leq K\|\eta\|\,\|\nabla U\|^2.
		\end{equation*}
	\end{lemma}
	
	\begin{proof}
		Notice that all boundary terms vanish since $\eta,\omega\in H_0^1(0,L)$. Then, integrating by parts gives
		\begin{align*}
			\int_0^L\eta(\eta\omega)_x\,dx 
			& = \frac12\int_0^L\omega_x\eta^2\,dx,\\
			\int_0^L\omega^2\omega_x\,dx&=0,\\
			\int_0^Lx\omega(\eta\omega)_x\,dx 
			& =-\int_0^L\eta\omega^2\,dx -\int_0^Lx\eta\omega\omega_x\,dx.
		\end{align*}
		Substituting $f=-\N(\eta,\omega)$ proves the identities. By Lemma \ref{lem:bil} and Poincar\'e's inequality, we reach
		\begin{align*}
			\frac12\left|\int_0^L\omega_x\eta^2\,dx\right|
			&\leq\frac12\|\omega_x\|\,\|\eta\|_{L^\infty}\|\eta\| \\
			&\leq\frac14\left(\frac L\pi\right)^{1/2}\|\eta\|\left(\|\eta_x\|^2+\|\omega_x\|^2\right),\\
			\left|\int_0^L\eta\omega^2\,dx\right|
			&\leq\|\eta\|\,\|\omega\|_{L^\infty}\|\omega\|\leq\left(\frac L\pi\right)^{3/2}\|\eta\|\,\|\omega_x\|^2.
		\end{align*}
		The estimate follows from the definition of $K$ in~\eqref{eq:cts}.
	\end{proof}
	
	Now, with the previous estimates in hand, we can prove our main Theorem. For simplicity, we are going to split the proof in three steps
	
	\begin{proof}[Proof of Theorem \ref{th:main}]
		Let $U$ be the maximal solution given by Theorem~\ref{th:local}.
		
		\medskip
		\noindent\textit{Step 1. Lyapunov estimate.}
		For every $T<T_{\max}$, Lemma~\ref{lem:bil} gives $-\N(\eta,\omega)\in L^1(0,T;X_0)$. Proposition~\ref{pr:lin} therefore yields $V\in W^{1,1}(0,T)$. Combining \eqref{eq:Id1}--\eqref{eq:Id3} with Lemma~\ref{lem:NL}, we obtain
		\begin{align*}
			V'+\lambda V
			&\leq \frac12
			\begin{pmatrix}a&b\end{pmatrix} \Psi \begin{pmatrix}a\\b\end{pmatrix}\\
			&\quad	-\frac{|\beta|}{2}(d-\dot\tau)b^2
			+\frac{\mu_1+\lambda}{2}\|U\|^2 + \lambda\mu_1V_1\\
			&\quad-\frac{3\mu_1}{2}\|\nabla U\|^2+K\|\eta\|\,\|\nabla U\|^2\\
			&\quad+\frac{|\beta|}{2}\int_0^1
			\Big[-\mu_2(1-\dot\tau\rho)\\
			&\qquad +\lambda\tau\big(1+\mu_2(1-\rho)\big)\Big]z^2\,d\rho,
		\end{align*}
		for almost every $t\in(0,T_{\max})$. The first two terms are nonpositive by Definition~\ref{def:adm} and \eqref{eq:tau}. Now, recalling that
		\begin{gather*}
			|V_1|\leq\frac L2\|U\|^2,
			\ \ 
			\|U\|^2\leq\frac{L^2}{\pi^2}\|\nabla U\|^2,\ \ 
			1-\dot\tau\rho\geq1-d,
		\end{gather*}
		together with $\tau\leq M$, we can claim that
		\begin{align*}
			V'+\lambda V
			&\leq
			\Big[\frac{L^2}{2\pi^2}\left(\mu_1+\lambda(1+\mu_1L)\right)  -\frac{3\mu_1}{2}\\
			&\qquad+K\|\eta\| \Big]\|\nabla U\|^2\\
			&\quad+\frac{|\beta|}{2}
			\left[\lambda M(1+\mu_2)-\mu_2(1-d)\right] \|z\|_{L^2(0,1)}^2.
		\end{align*}
		By \eqref{eq:lambda} and the definition of $\delta$ in \eqref{eq:cts}, we conclude that
		\begin{equation}\label{eq:key}
			V'(t)+\lambda V(t)
			\leq\left(-\delta+K\|\eta(t)\|\right)\|\nabla U(t)\|^2
		\end{equation}
		for almost every $t\in(0,T_{\max})$.
		
		\medskip
		\noindent\textit{Step 2. Energy decay and global existence.}
		The smallness assumption~\eqref{eq:smallness} gives
		\begin{equation*}
			\|\eta(0)\|\leq\sqrt{2E(0)}
			\leq\sqrt{2\zeta E(0)}<\frac{\delta}{K}.
		\end{equation*}
		As long as $\|\eta(t)\|<\delta/K$, integrating \eqref{eq:key} and using \eqref{eq:EeqV} yields
		\begin{equation}\label{eq:dec}
			E(t)\leq\frac{V(0)}{1-\mu_1L}e^{-\lambda t}
			\leq\zeta E(0)e^{-\lambda t}.
		\end{equation}
		If $t_*<T_{\max}$ were the first time at which $\|\eta(t_*)\|=\delta/K$, continuity and \eqref{eq:dec} implies
		\begin{equation*}
			\frac{\delta}{K} = \|\eta(t_*)\|
			\leq\sqrt{2\zeta E(0)}\,e^{-\lambda t_*/2} <\frac{\delta}{K},
		\end{equation*}
		a contradiction. Thus \eqref{eq:dec} holds on $[0,T_{\max})$. In particular,
		\begin{equation*}
			\sup_{0\leq t<T_{\max}}\|U(t)\|_t^2 = \sup_{0\leq t<T_{\max}}2E(t)
			\leq2\zeta E(0),
		\end{equation*}
		so Theorem~\ref{th:local}(c) gives $T_{\max}=\infty$.
		
		\medskip
		\noindent\textit{Step 3. Uniform bound.}
		Since $K\|\eta(t)\|\leq K\sqrt{2\zeta E(0)}<\delta$, integrating \eqref{eq:key} gives
		\begin{multline*}
			V(t)+\lambda\int_0^t V(r)\,dr\\
			+ \left(\delta-K\sqrt{2\zeta E(0)}\right) \int_0^t\|\nabla U(r)\|^2\,dr
			\leq V(0).
		\end{multline*}
		Using $V\geq0$ and \eqref{eq:EeqV}, and letting $t\to\infty$, we obtain
		\begin{equation*}
			\int_0^\infty\|\nabla U(t)\|^2\,dt
			\leq \frac{\left(1+\max\{\mu_1L,\mu_2\}\right)E(0)} {\delta-K\sqrt{2\zeta E(0)}}.
		\end{equation*}
		Together with the bound
		\begin{equation*}
			\sup_{t\geq0}\|U(t)\|\leq\sqrt{2\zeta E(0)},
		\end{equation*}
		this proves~\eqref{eq:mest}.
	\end{proof}
	
	\begin{remark}[A perturbative proof]\label{rem:pert}
		Local exponential stability can also be deduced from the linear decay by a perturbation argument, in the spirit of \cite{BCV,PTV}. We sketch it. Fix $T>0$ such that $\zeta e^{-\lambda T}\leq\frac14$, and let $U_0$ satisfy $$24C_T^2T^{1/4}\|U_0\|_0\leq\frac14.$$ By Theorem~\ref{th:local}(a), the solution $U$ exists on $[0,T]$ and $\|(\eta,\omega)\|_{B_T}\leq2C_T\|U_0\|_0$. Write $U=U_{\rm lin}+W$, where $U_{\rm lin}(t)=S(t,0)U_0$ solves the linearized system and $W$ is the mild solution of \eqref{eq:lin.sou} with zero initial state and source $-\N(\eta,\omega)$. Since $\|U_{\rm lin}(T)\|_T^2$ is twice the energy of $U_{\rm lin}$ at time $T$, \cite[Theorem 1.1]{CCGM}, with the corrections of Section~\ref{sec:corr}, gives
		\begin{equation*}
			\|U_{\rm lin}(T)\|_T\leq\left(\zeta e^{-\lambda T}\right)^{1/2}\|U_0\|_0\leq\frac12\|U_0\|_0,
		\end{equation*}
		while Proposition~\ref{pr:lin} and \eqref{eq:NL.est} give
		\begin{equation*}
			\|W(T)\|_T\leq 6T^{1/4}\|(\eta,\omega)\|_{B_T}^2\leq24C_T^2T^{1/4}\|U_0\|_0^2.
		\end{equation*}
		Hence
		\begin{equation*}
			\|U(T)\|_T \leq\left(\frac12+24C_T^2T^{1/4}\|U_0\|_0\right)\|U_0\|_0\leq\frac34\|U_0\|_0.
		\end{equation*}
		In particular, $U(T)$ satisfies the same smallness condition as $U_0$ and, by Lemma~\ref{lem:shift}, the problem can be restarted at time $T$ with the same constants. Iterating, $\|U(nT)\|_{nT}\leq(3/4)^n\|U_0\|_0$ for every $n\geq1$, which yields the exponential decay of $E$. This argument provides a decay rate depending on $T$ and a less explicit smallness condition, whereas the Lyapunov proof retains the linear rate $\lambda$ and the explicit threshold \eqref{eq:smallness}.
	\end{remark}
	
	\section{Numerics and concluding remarks}\label{sec:conclusions}
	
	\subsection{Numerical illustration and verification}
	
	We illustrate Theorem~\ref{th:main} for system~\eqref{eq:KdV} with $L=5$, $\alpha=1$, $\beta=1/2$ and $\tau(t)=2+\sin(t/2)$, so that $1\leq\tau(t)\leq3$ and $\dot\tau(t)\leq1/2$. The admissible choice
	\begin{equation*}
		(\mu_1,\mu_2)\simeq(0.04777231,\,0.04455379),
	\end{equation*}
	with $\lambda=\lambda^*/2\simeq0.00355445$, gives $\zeta\simeq1.62764$ and $r^2\simeq5.64783\times10^{-5}$. For zero initial history, we take
	\begin{equation*}
		\eta_0(x)=A s^2(1-s)^2(1+2s),\qquad
		\omega_0(x)=A s(1-s)^2,
	\end{equation*}
	with $s=x/L$,
	where $A>0$ is chosen so that $E(0)=\kappa r^2$ for $\kappa\in\{0.1,0.5,0.9\}$. These smooth data satisfy the boundary compatibility conditions and the strict smallness assumption $E(0)<r^2$.
	
	We use Legendre--Galerkin spaces of degrees $N$ and $N-1$ for $\eta$ and $\omega$, respectively, imposing the  boundary conditions and the feedback through the boundary term. Gauss quadrature integrates the polynomial terms exactly. Time integration uses second-order semi-implicit backward differentiation, with explicit extrapolation of the nonlinearity and a backward-Euler starting step. The delayed trace is linearly interpolated from the prescribed delayed profile and previously computed values. The primary runs use $N=24$, $\Delta t=0.005$ and $T=400$. Writing $a_N(t)=\eta_{N,x}(t,L)$, we compute
	\begin{multline*}
		E_N(t)=\frac12\int_0^L(\eta_N^2+\omega_N^2)\,dx+ \frac{|\beta|}{2}\int_{t-\tau(t)}^t a_N(s)^2\,ds.
	\end{multline*}
	
\begin{figure*}[t]
		\centering
		\includegraphics[width=\textwidth]{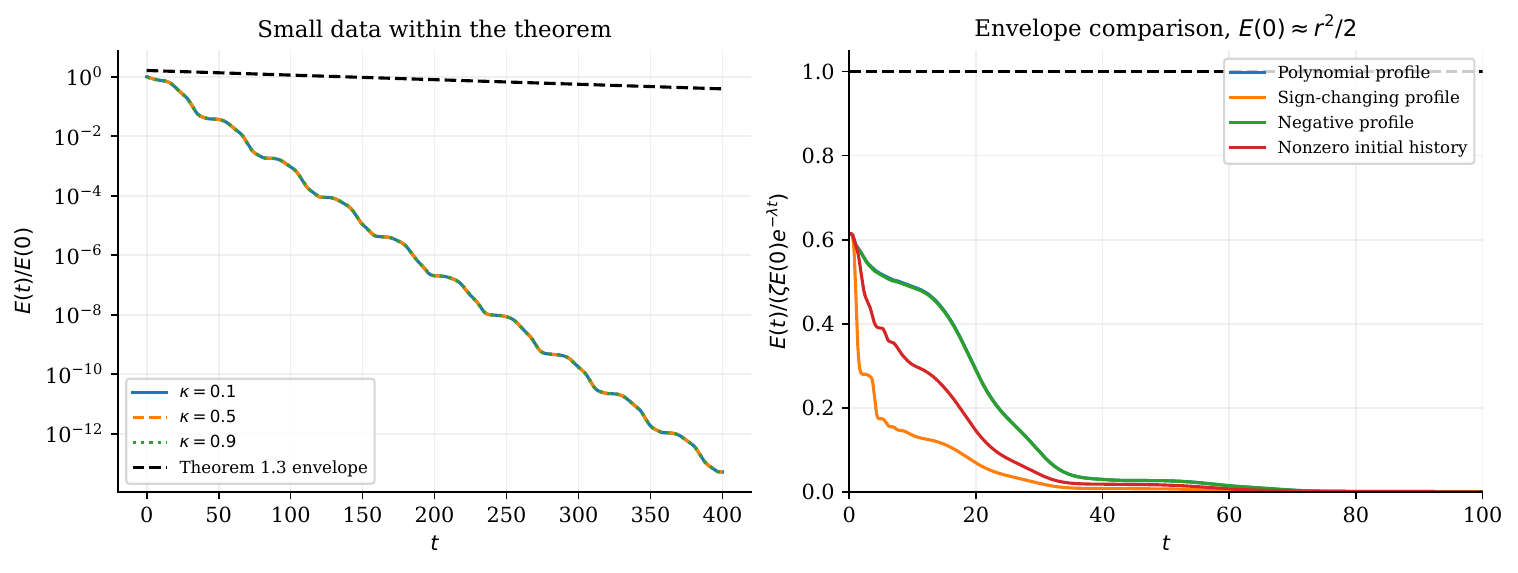}
		\caption{Small-data experiments with $N=24$, $\Delta t=0.005$ and  $\tau(t)=2+\sin(t/2)$. Left: normalized energies for  $\kappa=0.1,0.5,0.9$ and the theoretical envelope  $\zeta e^{-\lambda t}$. Right: ratios to the envelope for four  compatible initial states with $E_N(0)\simeq r^2/2$, shown on $[0,100]$.
			The horizontal dashed line marks equality with the envelope.}
		\label{fig:num-small}
	\end{figure*}

	Figure~\ref{fig:num-small} compares the computed energies with $\zeta E_N(0)e^{-\lambda t}$. The normalized trajectories almost coincide, and their maximum ratios to this envelope are below $0.615$ on $[0,400]$. Additional experiments with the negative primary profile, a sign-changing profile and compatible nonzero history, at total energy approximately $r^2/2$, also remain below the envelope; their definitions are supplied in \cite{kdv-delay-notebook}. 
	
\begin{figure*}[t]
		\centering
		\includegraphics[width=\textwidth]{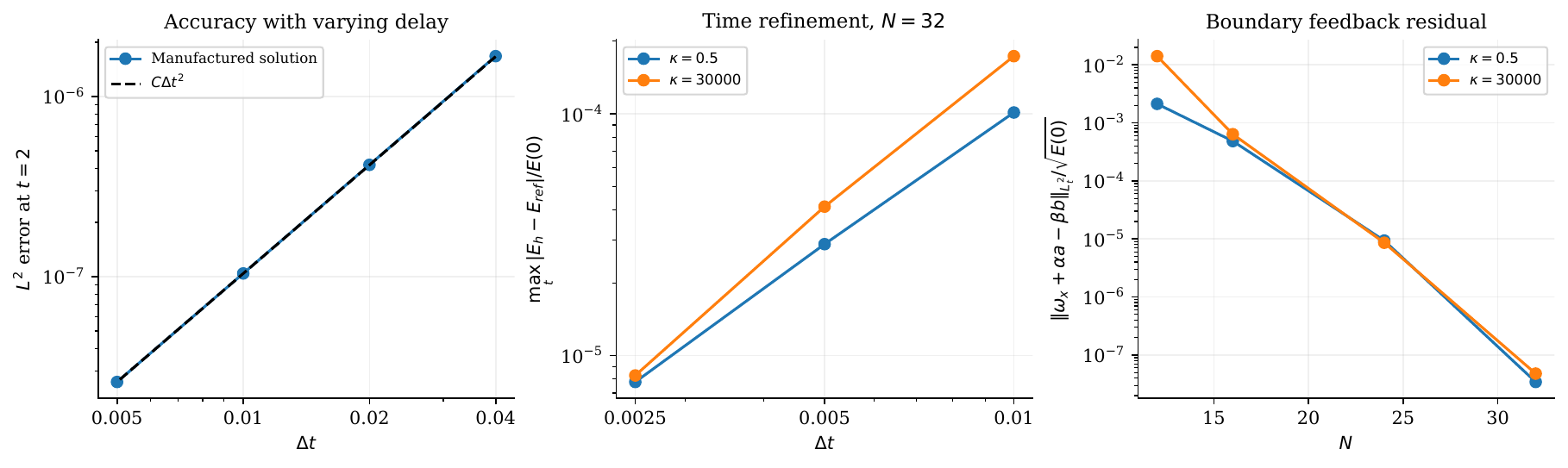}
		\caption{Numerical verification. Left: manufactured-solution $L^2$  errors at $t=2$ and a second-order reference line. Centre: energy  differences on $[0,100]$ at $N=32$, relative to $\Delta t=0.00125$.  Right: normalized time-$L^2$ boundary-feedback residuals at  $\Delta t=0.0025$. 
			The refinement curves correspond to  $\kappa=0.5$ and $30000$; the latter is a test of numerical resolution outside the theorem's small-data regime.}
		\label{fig:num-validation}
	\end{figure*}

Figure~\ref{fig:num-validation} reports numerical resolution checks.
A manufactured nonlinear solution with
a nonzero delayed profile and time-varying
delay yields temporal orders between
$2.001$ and $2.005$.. A manufactured nonlinear solution with a nonzero delayed profile and time-varying delay yields temporal orders between $2.001$ and $2.005$.
	For $\kappa=0.5$, the maximum energy difference on $[0,100]$, normalized by $E_N(0)$, is $2.89\times10^{-5}$ when comparing $\Delta t=0.005$ with $0.00125$ at $N=32$, and $9.21\times10^{-10}$ when comparing $N=24$ with $32$ at $\Delta t=0.0025$. The boundary-feedback residual also decreases with $N$. The larger-amplitude curves test numerical resolution outside the theorem's range. These findings support the reported computations; they are neither rigorous error bounds, they are mainly designed to verify the findings without being a numerical test.
	
	\subsection{Final comments}
	\begin{enumerate}
		\item
		Theorem~\ref{th:main} establishes global existence and exponential decay for sufficiently small initial data, with the same admissible range of decay rates as in the linear estimate. Global existence for arbitrary data remains open for the feedback considered here. The difficulty lies in the cubic contribution $-\frac12\int_0^L\omega_x\eta^2\,dx$, which has no prescribed sign and prevents the present estimates from providing a global energy bound without a smallness assumption. This difficulty persists in the absence of delay.
		
		\item
		The decay rate and the smallness threshold depend on the admissible pair $(\mu_1,\mu_2)$. For a fixed pair, increasing $\lambda$ decreases $\delta$ and therefore reduces the threshold in~\eqref{eq:smallness}. Maximizing the decay bound and maximizing the size of the guaranteed neighborhood of attraction are thus different optimization problems. As shown in (C8) of Section~\ref{sec:corr}, optimization along the curve of (C4) need not give the best decay bound over the full admissible set. These Lyapunov estimates do not determine the optimal decay rate of the system or the full basin of attraction.
		
		\item
		The restriction $L<\sqrt3\pi$ arises from Poincar\'e's inequality in the Lyapunov estimate. The stability result cannot hold for every length: if $L=2\pi n$, with $n\geq1$, then
		$(\eta,\omega,z)=\left(c(1-\cos x),0,0\right)$,  $c\in\mathbb R\setminus\{0\}$,
		is a stationary solution of~\eqref{eq:KdV}, with zero initial delayed profile, for every delay satisfying~\eqref{eq:tau}. Indeed, $\eta_x+\eta_{xxx}=0$ and $\eta=\eta_x=0$ at both endpoints. These states have energy $3c^2L/4$, which can be arbitrarily small, so the zero state is not locally attractive at these lengths. Determining which other lengths admit small-data exponential stability remains a further question.
	\end{enumerate}

	\section*{Declaration of competing interest}
	The authors declare that they have no known competing financial interests or personal relationships that could have appeared to influence the work reported in this paper.
	
	\section*{Funding}
	The Article Processing Charge (APC) for the publication of this research was funded by the Coordenação de Aperfeiçoamento de Pessoal de Nível Superior -Brasil (CAPES) (ROR identifier: 00x0ma614).

 This work was completed while J.~R.~Mu\~noz was visiting R.~de~A.~Capistrano--Filho and V.~H.~Gonzalez~Martinez at UFPE. J.~R.~Mu\~noz acknowledges the support of ANID, Chile, through FONDECYT Postdoctorado Grant No.~3260141. R.~de~A.~Capistrano--Filho was partially supported by CAPES/CO\-FE\-CUB grant number 88887.879175/2023-00, CNPq grant numbers 301744/2025-4 and 421573/2023-6, and PROPG (UFPE) \textit{via} PROAP resources. V.~H.~Gonzalez~Martinez was supported by CAPES/COFECUB grant 88887.879175/2023-00 and CNPq grants 421573/2023-6, 402652/2024-0 and 300997/2025-6.

	\section*{CRediT authorship contribution statement}
	\textbf{Roberto de A. Capistrano--Filho:} Conceptualization, Formal analysis, Writing -- original draft. \textbf{Boumedi\`ene Chentouf:} Formal analysis, Writing -- review \& editing. \textbf{Victor H. Gonzalez~Martinez:} Formal analysis, Writing -- review \& editing. \textbf{Juan Ricardo Mu\~noz:} Formal analysis, Software, Validation, Writing -- review \& editing.
	
	\section*{Data availability}
	The code used for the numerical illustrations in Section~\ref{sec:conclusions} is available in \cite{kdv-delay-notebook}.

\end{document}